\documentclass[11pt,letterpaper]{amsart}
\usepackage[english]{babel}
\usepackage{amssymb}
\usepackage{fourier}
\usepackage{mathrsfs}
\usepackage{enumerate}
\usepackage{color}
\usepackage{ifpdf}
\usepackage{tikz}
\usepackage{tikz-cd}
\usetikzlibrary{decorations.pathmorphing,arrows}
\usepackage[initials]{amsrefs}

\usepackage{caption}
\usepackage{subcaption}
\usepackage{comment}
\usetikzlibrary{intersections}

\newcommand{\bbN}{{\mathbb N}}
\newcommand{\bbQ}{{\mathbb Q}}
\newcommand{\bbR}{{\mathbb R}}

\newcommand{\bbZ}{{\mathbb Z}}

\newcommand{\bfH}{{\mathbf H}}

\newcommand{\GL}{\operatorname{GL}}

\newcommand{\SL}{\operatorname{SL}}
\newcommand{\PSL}{\operatorname{PSL}}

\newcommand{\GCD}{\mathrm{GCD}}
\newcommand{\LCM}{\mathrm{LCM}}

\newcommand{\setdef}[2]{ \left\{ {#1} \left. \mid {#2} \right.\right\} }

\newtheorem{theorem}{Theorem}[section]
\newtheorem{thm}[theorem]{Theorem}
\newtheorem{mthm}{Theorem}

\newtheorem{lemma}[theorem]{Lemma}

\newtheorem{corollary}[theorem]{Corollary}

\newtheorem{conjecture}[theorem]{Conjecture}

\newtheorem{claim}[theorem]{Claim}
\newtheorem{problem}[theorem]{Problem}

\theoremstyle{definition}
\newtheorem{definition}[theorem]{Definition}

\newtheorem{example}[theorem]{Example}

\newtheorem{remark}[theorem]{Remark}

\numberwithin{equation}{section}

\title[Trace set of hyperbolic surfaces]{On trace sets of hyperbolic surfaces and a conjecture of Sarnak and Schmutz}
\author{Yanlong Hao}
\address{University of Michigan, Ann Arbor, MI, USA}
\email{ylhao@umich.edu}

\begin{document}
\date{\today}

\begin{abstract}
In this paper, we investigate the trace set of a Fuchsian lattice. There are two results of this paper: The first is that for a non-uniform lattice, we prove Schmutz’s conjecture: the trace set of a Fuchsian lattice exhibits linear growth if and only if the lattice is arithmetic. Additionally, we show that for a fixed surface group $\Gamma_g$ of genus $g\geq 3$ and any $\epsilon>0$, the set of cocompact lattice embeddings such that their growth rate of the trace set exceeds $n^{2-\epsilon}$ has positive Weil-Petersson volume. We also provide an asymptotic analysis of the volume of this set as $g\to \infty$.
\end{abstract}
\maketitle

\section{Introduction and Statement of the main results} 
We say that a set $A$ of real numbers satisfies the \emph{bounded clustering}  \emph{(B-C) property}
if and only if  there exists a constant $K_A$ such that $A\cap [n,n+1]$ has less than $K_A$ elements for all $n\in \bbZ$. Furthermore, set
$${\text{Gap}}(A):=\inf\setdef{|a-b|}{a,b\in A, a\neq b}$$

Let $\Gamma$ be a subgroup of $\PSL(2,\mathbb{R})$. The \emph{trace set} $\text{Tr}(\Gamma)$ of $\Gamma$ is defined (up to a sign) as the set of traces of elements of $\Gamma$. 

In \cite{luo1994number}, Luo and Sarnak showed that the trace set of an arithmetic Fuchsian group satisfies the B-C property. Recalled that a Fuchsian group is a discrete subgroup of $\PSL(2,\mathbb{R})$. Furthermore, Sarnak conjectured that the converse is also true.
\begin{conjecture}[Sarnak \cite{MR1321639}]\label{BC}
   Let $\Gamma$ be a cofinite Fuchsian group. 
   \begin{enumerate}
       \item 
   If $\mathrm{Tr}(\Gamma)$ satisfies the B-C property, then $\Gamma$ is arithmetic. 
   \item
   If $\mathrm{Gap(Tr}(\Gamma))>0$, then $\Gamma$ is derived from a quaternion algebra.
   \end{enumerate}
\end{conjecture}

In \cite{MR1394753}, Schmutz makes an even stronger conjecture using the linear growth of a set instead of the B-C property. A subset of reals is said to have linear growth if and only if there exist positive real constants $C$ and $D$ such that for all $n$,
$$\#\setdef{a\in A}{|a|\leq n}\leq Cn+D.$$

\begin{conjecture}[Schmutz \cite{MR1394753}]\label{LG}
    Let $\Gamma$ be a cofinite Fuchsian group. If
${\mathrm{Tr}}(\Gamma)$ has linear growth then $\Gamma$ is arithmetic. 
\end{conjecture}
 Schmutz proposed a proof of Conjecture~\ref{LG} for nonuniform lattices. In \cite{MR1394753}, Schmutz essentially proved part (2) of Sarnak's Conjecture~\ref{BC} under (1). Unfortunately, the proof of Conjecture~\ref{LG} contains a gap. Later, Geninska and Leuzinger fixed part of this gap in \cite{geninska2008geometric} and confirmed part (1) of Sarnak's Conjecture for nonuniform Fuchsian lattices.  Note that Conjecture~\ref{LG} is still open even for nonuniform lattices. And the Conjecture~\ref{BC} remains entirely open for cocompact Fuchsian groups.

In this paper, we first prove a stronger version of Conjecture~\ref{LG} for non-uniform lattices. 
\begin{mthm}\label{main}
Let $\Gamma$ be a non-uniform lattice of  $\PSL(2,\mathbb{R})$. 
   If 
   $$\lim_{n\to \infty}\frac{\#\{\mathrm{Tr}(\Gamma)\cap [-n,n]\}}{n\log\log\log n}=0,$$ then $\Gamma$ is arithmetic. 
\end{mthm}

There is also a geometrical version of Theorem~\ref{main}. Let $\Gamma$ be a torsion-free non-uniform lattice of  $\PSL(2,\mathbb{R})$. Then ${\bf H}^2/\Gamma$ is a complete hyperbolic surface, where ${\bf H}^2$ is the upper-half plane. Each hyperbolic conjugacy class 
$[\gamma]$ corresponds to a unique closed geodesic on ${\bf H}^2/\Gamma$. We denote the length of the geodesic related to $[\gamma]$ by $\ell(\gamma)$. Then we have 
\[|\mathrm{Tr}(\gamma)|=2\cosh(\frac{\ell(\gamma)}{2}).\]

For any complete hyperbolic surface $(\Sigma, d)$, we denote the length set $LS_d$ of $(\Sigma,d)$ to be the set of all lengths of closed geodesics on $(\Sigma,d)$ without multiplicities. A direct translation of Theorem~\ref{main} has the following form. 
\begin{thm}\label{LS}
    Let $(\Sigma,d)$ be a complete, non-compact hyperbolic surface of finite volume. 
    If 
    \[\lim_{n\to \infty}\frac{\#(LS_d\cap[0,n])}{e^{\frac{n}{2}}\log\log n}=0,\]
then $(\Sigma,d)$ is arithmetic.
\end{thm}

Theorem~\ref{main} is a consequence of the following result.
\begin{mthm}\label{main2}
    Let $\Gamma$ be a subgroup of $SL_2(\mathbb{R})$. Assume that $\Gamma$ contains a parabolic element, and 
    $$\lim_{n\to \infty}\frac{\#\{\mathrm{Tr}(\Gamma)\cap [-n,n]\}}{n\log\log\log n}=0.$$
    Then up to conjugation, one of the following is true:
    \begin{enumerate}
        \item $\Gamma$ is elementary. Then, $\Gamma$ is upper-triangular, and the projection to the upper-left entry is
a homomorphism $\phi_{11}: \Gamma\to (\bbR^*, \cdot)$. The image of $\phi_{11}$, $D(\Gamma)$ is one of the following groups: the trivial group, $\bbZ/2\bbZ$, $\bbZ$, or $\bbZ\times \bbZ/2\bbZ$. If $D(\Gamma)$ has a $\bbZ$-factor, let $\lambda$ be a generator of that factor, and set $\gamma=\begin{pmatrix}
    \lambda& 0\\
    0&\frac{1}{\lambda}
\end{pmatrix}\in\Gamma.$ If $D(\Gamma)$ has no $\bbZ$-factor, let $\lambda=1$ and $\gamma=\begin{pmatrix}
    1& 0\\
    0&1
\end{pmatrix}$. Then $\Gamma$ has the form $A\ltimes Q$, where $A=\langle \gamma\rangle$ and $Q=\phi_{11}^{-1}\{\pm 1\}$. Furthermore, let $P=\phi_{11}^{-1}\{1\}$. Then $Q$ has the following forms: 
\begin{enumerate}
    \item If $D(\Gamma)$ has no $\bbZ/2\bbZ$-factor, then $Q=P:=\{\begin{pmatrix}
       1&t\\
       0&1
    \end{pmatrix}\mid t\in R\}$ for a (possibly trivial) subgroup $R$ of $\mathbb{R}$ satisfying $\lambda^2R=R$.
    \item\label{b}  If $D(\Gamma)$ has a $\bbZ/2\bbZ$-factor, then there is a short exact sequence
    \[
    1\to P\to Q\xrightarrow{\phi_{11}} \{\pm 1\}\cong\bbZ/2\bbZ\to 1,
    \]
    where $P$ has the form described above, and the set of upper-right entries of
$Q$ is also invariant under multiplication by $\lambda^2$.
\end{enumerate}

        \item $\Gamma$ is non-elementary. Then $\Gamma$ is discrete, and there is a finite index subgroup $\Gamma'$ of $\Gamma$ which conjugates to a subgroup of $\SL(2,\mathbb{Z}).$
    \end{enumerate}
\end{mthm}
\begin{example}
    Let $\Gamma=\{\begin{pmatrix}
        -1& k+\frac{1}{2}\\
        0&-1
    \end{pmatrix}\mid k\in\bbZ\}\bigcup \{\begin{pmatrix}
        1& k\\
        0& 1
    \end{pmatrix}\mid k\in\bbZ\}$. Then $\Gamma$ is an elementary Fuchsian group and belongs to case~\ref{b}. Note that $\Gamma$ contains no subgroup isomorphic to $\bbZ/2\bbZ$.
\end{example}\label{example: 1b}
By considering the growth function of cases in Theorem~\ref{main2}, we have the following.
\begin{corollary}
    Let $\Gamma$ be a subgroup of $SL_2(\mathbb{R})$. Assume that $\Gamma$ contains a parabolic element, and 
    $$\lim_{n\to \infty}\frac{\#\{\mathrm{Tr}(\Gamma)\cap [-n,n]\}}{n\log\log\log n}=0.$$ 
    Denote $f_\Gamma(n)=\#\{\mathrm{Tr}(\Gamma)\cap [-n,n]\}$, then one of the following is true:
    \[f_\Gamma(n)\sim 1, \quad or \quad f_{\Gamma}(n)\sim \log n, \quad or \quad f_{\Gamma}(n)\sim n.\]
\end{corollary}
Hence, there are gaps in the growth rate of the trace sets of subgroups of $\PSL(2,\mathbb{R})$ with parabolic elements.

\medskip

Now we turn to the proof of Theorem~\ref{main2}. The proof follows from an idea of Schmutz.

 First, we recall Schmutz's idea to prove  Conjecture~\ref{BC} here. Let $\Gamma$ be a non-uniform (torsion-free) Fuchsian lattice. Given an element $\gamma$ in $\Gamma$, Schmutz constructs a $Y$-piece $S$, a surface of signature (0,3) related to $\gamma$. By considering traces of different families of elements in $\pi_1(S)$, there are restrictions on the trace of $\gamma$. This implies that the lattice commensurates to a subgroup of $\PSL(2,\mathbb{Q})$. After this stage, the construction by Geninska and Leuzinger continues the argument, restricting the trace set to be a subset of $\mathbb{Z}$ under the B-C property, and completes the proof.

Now, when $\Gamma$ is elementary, the classification in Theorem~\ref{main2} follows from an easy argument.

In the non-elementary case, the proof of Theorem~\ref{main2} has two steps. We follow a strategy similar to Schmutz’s proposal, but our approach is more algebraic. It is essentially the translation of Schmutz's work into algebraic language. But we still state it in full since in Theorem~\ref{main2}, we deal with general subgroups and lack of geometrical picture. Under weaker assumption as in Theorem~\ref{main2}, we can restricted the group $\Gamma$ to a related normal group ${\Gamma^{(2)}}$  which is a subgroup of $\PSL(2,\mathbb Q)$. Indeed, we would like to work with a slightly bigger but still normal subgroup $\bar{\Gamma}$. 

We then replace the work of Geninska and Leuzinger with a new approach: note that through left/right multiplication by unipotents, each non-trivial element $\gamma\in\bar{\Gamma}$ generates a subset $A$ of traces with non-trivial density, i.e.
$$\lim_{n\to\infty}\frac{\#\{A\cap [-n,n]\}}{2n}>0.$$
If there is an element in $\Gamma$ whose trace belongs to $\mathbb{Q}$ but not in $\mathbb{Z}$, we construct a countable family of subsets of the trace set such that their union has infinite density. The results follow by considering the intersection of this union with $[-n,n]$. It follows that $\mathrm{Tr}(\bar{\Gamma})\subset \mathbb{Z}$. Since $\bar{\Gamma}$ contains two linearly independent unipotent elements, $\bar{\Gamma}$ is, up to conjugation and passing to a finite index subgroup, a subgroup of $\SL(2,\mathbb{Z})$. The final step is to show that the quotient group $\Gamma/\bar{\Gamma}$ is finite.

\medskip

We now investigate Conjecture~\ref{BC} and $\ref{LG}$ for cocompact Fuchsian lattices. Although we cannot completely resolve these conjectures, we assert that, for large genus $g$, the conjectures hold for a big subset of the Teichm\"uller space $\mathcal{T}_g$.

First, we present a method for constructing hyperbolic surfaces with a large trace set. Let $(\gamma_i)_{i=1}^{2g}$ be a standard generating set of $\Gamma_g$. The free subgroup generated by $(\gamma_i)_{i=1}^{2g-3}$ is called a $(2g-3)$-subgroup if it is a convex-cocompact Fuchsian subgroup. 

For  a fixed $(2g-3)$-subgroup $A$, let $\delta_A: \mathcal{T}_g\to \mathbb{R}^+$ such that $\delta_A([d])$ is the critical exponent of the Fuchsian group $\psi_d(A)$ where $\psi_d$ is the lattice embedding corresponding to $[d]$.

\begin{mthm}\label{main: cocompact 2}
    Let $\Sigma_g$ be a closed surface of genus $g\geq 3$, and let $\Gamma_g$ 
 be its fundamental group. Let $A$ be a $(2g-3)$-subgroup of $\Gamma_g$. For any point $[d]\in\mathcal{T}_g$ and any $\epsilon>0$, there exists a neighborhood $V_g^\epsilon(d)$ of $[d]$ and a subset $T_{\mathrm{sing}}\subset \mathcal{T}_g$, which is a union of countably many algebraic subsets of positive codimension in $\mathcal{T}_g$ and is independent of $\epsilon$, such that for all hyperbolic metric $[d']\in V_g^\epsilon(d)\setminus T_{\mathrm{sing}}$, the corresponding lattice has trace growth greater than $n^{2\delta_{A}(A)-\epsilon}$.
\end{mthm}

\begin{remark}
    The set $T_{\mathrm{sing}}$ is exactly the same set as the $T_{\mathrm{sing}}$ in \cite{hao2022marked}*{Theorem D}.
\end{remark}

The strategy of proving Theorem~\ref{main: cocompact 2} is based on the belief that, for general points in the Teichmuller space, the multiplicities of their trace sets are minimal among all complete hyperbolic structures. Therefore, if we can find one hyperbolic structure where the growth rate of the trace set is large, then, in an open neighborhood of that point, almost all points will have a large growth rate of their trace sets.

Thus, to find a lattice with a trace growth rate greater than $n$, it is sufficient to identify a point $[d]\in \mathcal{T}_g$ such that $\delta_{\psi_d(A)}>\frac{1}{2}$. This is achieved by selecting short separating multicurves on the surface and considering the geometric Cheeger constant.

Let $\Sigma_g$ be a closed surface of genus $g\geq 3$, and let $\Gamma_g$ be its fundamental group. Since the critical exponent of $A$ is equivariant under the action of the mapping class group (the subgroup $A$ changes according to different generating sets). We can consider the Moduli space for $\Sigma_g$, denoted by $\mathcal{M}_g$. 

Following \cite{mirzakhani2013growth}, we consider the following objects: Let $\Xi_2(g)$ be the set of multicurves $\alpha$ on $\Sigma_g$, where $\alpha=\cup_{i=1}^s\alpha_i$, such that all $\alpha_i$ are simple closed geodesics, and $\Sigma_g\setminus \alpha=\Sigma_1\cup \Sigma_2$, with $\Sigma_1$ and $\Sigma_2$ being connected subsurfaces and $|\chi(\Sigma_2)|=2\leq |\chi(\Sigma_1)|.$ Here, $|\chi(\Sigma_1)|=2g_1-2+s$ is the absolute value of the Euler characteristic. 

Computations show that  $|\chi(\Sigma_1)|=2g-4$, meaning that $\pi_1(\Sigma_1)$ is a free group of rank $2g-3$. On the other hand, $\Sigma_2$ has surface type $(0,4)$ or $(1,2)$. A detailed verification shows that there exists a standard generating set of $\Gamma_g$ such that the first $2g-3$ elements generate the fundamental group of $\Sigma_1$. We take $A:=\pi_1(\Sigma_1)$ as a subgroup of $\Gamma_g$ induced by the surface embedding. We call such a subgroup $A$ \textit{geometrically selected.} Note that $A$ is a $(2g-3)$-group.

When the length of the multicurve $\alpha$ is short, by considering the Cheeger constant for the surface ${\bfH}^2/A$, we show that the critical exponent of $A$ is large. 

Applying \cite{mirzakhani2013growth}*{Theorem 4.9}, we have the following:
\begin{mthm}\label{main: cocompact}
Consider the Weil-Petersson volume on the Moduli space $\mathcal{M}_g$. For any $\epsilon>0$, let $\bar{V_g^\epsilon}$ denote the set of points in $\mathcal{M}_g$ such that the corresponding lattice has trace growth greater than $n^{2-\epsilon}$. Then 
\[\lim_{g\to \infty} \frac{\mathrm{Vol_{WP}}(\bar{V_g}^{\frac{\epsilon}{g}})}{\mathrm{Vol_{WP}}(\mathcal{M}_g)}=1.\]
\end{mthm}

Let $(\Sigma_g, d)$ be a negatively curved Riemannian surface and $\Gamma_g$ its fundamental group. We normalize the metric so that the topological entropy of $d$ is 1. We define the trace set of $d$ as 
\[\mathrm{Tr}(d):=\{2\cosh(\frac{\ell_d({\gamma})}{2})\mid \gamma\in\Gamma\},\] where $\ell_d$ is the marked length spectrum, which gives the length of the closed geodesic representing the conjugacy class $[\gamma]$. Since the length spectrum of a generic negatively curved Riemannian manifold has multiplicity 1 \cite{abraham1970bumpy}, by Margulis' prime geodesic theorem for variance curvature manifolds \cite{Margulis-thesis}, we have that for a generic negatively curved metric $d$ on $\Sigma_g$:
$$\lim_{n\to \infty} \frac{\#\{\mathrm{Tr}(d)\cap [0,n]\}}{\frac{n^2}{\log n}}>0.$$
However, the multiplicities of the length set are unbounded for hyperbolic metrics, \cite{randol1980length}, and the asymptotic (average) of these multiplicities remains largely unknown (see some progress in \cites{ginzburg1998stable, leininger2003equivalent}). Note that Buser in \cite{buser1992geometry}*{Remark 3.7.13} has shown that there exists $c>0$ and a
sequence $l_n\in L(S)$ such that the multiplicity of $l_n$ is at least $cl_n^{\frac{\log 2}{\log 5}}$. Hence, the results of Theorem~\ref{main: cocompact 2} and ~\ref{main: cocompact} do not follow from earlier results.

\subsection{Related works} There are also related works on the structure or properties of the length set. In \cite{lafont2019primitive}, Lafont and McReynolds showed that every noncompact, locally symmetric, arithmetic manifold has arbitrarily long arithmetic progressions in its primitive length spectrum. This result was extended by Miller \cite{miller2016arithmetic} to every arithmetic locally symmetric orbifold of classical type without Euclidean or compact factors. A deeper structure for the length set of subarithmetic hyperbolic cusped manifolds is revealed in \cite{kontorovich2024length}.

Another direction of research concerns the rigidity problem. In \cite{chinburg2008geodesics} and \cite{prasad2009weakly}, the authors showed that length-commensurability has strong implications, one of which is that length-commensurable, arithmetically defined, locally symmetric spaces of certain types are necessarily commensurable.

There are also numerical results \cite{bogomolny2004multiplicities} suggest that the average multiplicities of the length spectrum have exponential growth (with a smaller exponent) for certain
non-arithmetic surfaces associated with Hecke triangle groups. Note that here the corresponding manifolds (orbifolds) are non-compact. 

\subsection{Further research and open questions}
Theorem~\ref{main2} raises similar questions for $\SL(2,\mathbb{C})$, which we will address in a forthcoming paper. 
\begin{problem}
    Generalise Theorem~\ref{main} and \ref{main2} for the trace sets of Kleinian groups, or higher-dimensional Kleinian groups with a suitable notation of trace.
\end{problem}

Moreover, the trace set contains more information than the length set of a higher-dimensional hyperbolic manifold. Hence, although generalizations of  Theorems~\ref{main} and \ref{main2} to 3-dimensional hyperbolic manifolds are likely, Theorem~\ref{LS} still requires further work.

\begin{problem}
    Let $(M,d)$ be a complete, non-compact hyperbolic 3-manifold of finite volume. 
    If 
    \[0<\liminf_{n\to \infty}\frac{\#(LS_d\cap[0,n])}{e^{n}}<\infty,\]
is $(M,d)$ arithmetic?
\end{problem}

For cocompact lattices, we propose the following conjecture:
\begin{conjecture}
 With the same notation as in Theorem~\ref{main: cocompact}, there exists a $0<\epsilon<1$ so that for sufficiently large genus $g$,
 \[\  \frac{\mathrm{Vol_{WP}}(\bar{V_g^\epsilon})}{\mathrm{Vol_{WP}}(\mathcal{M}_g)}=1.\]
\end{conjecture}

\medskip

The note is organized as follows: Section 2 recalls some preliminaries. In Section 3, we discuss the properties of second order linear recurrence sequences. Section 4 contains the proof of Theorem~\ref{main2} for Fuchsian groups, while Section 5 presents the proof of Theorems~\ref{main: cocompact 2} and ~\ref{main: cocompact}.

\medskip

\textbf{Acknowledgements.} I would like to acknowledge and thank Ralf Spatizer for his many suggestions, Beibei Liu for the discussion on the Cheeger constant for non-compact manifolds, and Peter Sarnak for providing reference \cite{bogomolny2004multiplicities}. I also want to thank the referees for their many helpful suggestions, and in particular for improving and correcting the statement of Theorem~\ref{main2}, for drawing attention to both Example~\ref{example: 1b} and a gap in the estimation (II) in Section~\ref{Step 2} of an earlier version.

\section{Definitions, notations, and some preliminaries}
First, we set up some basic notations used throughout this paper. For two functions $f$, $g:\mathbb{N}\to \mathbb{R}$, we say $f=O(g)$ if there exists $M>0$ such that $|f(n)|\leq M|g(n)|$ for all $n\in \mathbb{N}$. Additionally, $f=o(g)$ if $\lim_{n\to\infty}\frac{f(n)}{g(n)}=0$. Finally, $f\sim g$ if $f=O(g)$ and $g=O(f)$. 

We adopt the convention that if a function $f$ is only partially defined on $n>K$ for some $K$, we extend it to a function where $f(n)=1$ if $n$ is outside the domain, and continue to call it $f$ by abuse of notation. A typical example would be $f(n)=\log\log\log n$.
\subsection{Fuchsian groups}\label{Fuchsian} A general reference for this section is the book \cite{maclachlan2003arithmetic}. We denote by $\SL(2,\mathbb{R})$
the group of real $2\times2$ matrices with determinant 1, and by $\PSL(2,\mathbb{R})$ the
quotient group $\SL(2,\mathbb{R})/\{\pm I_2\}$ where $I_2$ is the $2\times 2$ identity matrix. 

A discrete subgroup of $\PSL(2,\mathbb{R})$ is called a \emph{Fuchsian group}. Let $p_{\mathbb R}:\SL(2,\mathbb{R})\rightarrow\PSL(2,\mathbb R)$ be the projection. Given a Fuchsian group $\Gamma$, we define the \emph{trace set} of $\Gamma$ as  
$$\mathrm{Tr}(\Gamma):=\setdef{\mathrm{tr}T}{T\in p_\mathbb{R}^{-1}{\Gamma}}.$$ 

A lattice of a locally compact, second countable topological group $G$ is a discrete subgroup $\Gamma$ such that $G/\Gamma$ has finite Haar measure. A lattice is called uniform if $G/\Gamma$ is compact, and nonuniform otherwise. 

A Fuchsian lattice $\Gamma$ is nonuniform if and only if $\Gamma$ contains parabolic elements. 

A Fuchsian group $\Gamma$ is called \textit{elementary} if it has a global fixed point in the compactification ${\bfH}^2\cup \partial {\bfH}^2$. Note that our definition differs somewhat from the standard one, since we regard focal subgroups as elementary. This broader viewpoint is appropriate for the aims of our work.

\medskip 

Arithmetic Fuchsian groups are obtained in the following way (see \cite{maclachlan1987commensurability} for example): Let $k$ be a totally real algebraic number field
with exactly one real archimedean place so that the $\bbQ$-isomorphisms of $k$ into
$\mathbb{C}$ are $\phi_1$, $\phi_2$, $\ldots$, $\phi_n$ where we take $\phi_1=\mathrm{Id}$, and $\phi_i(k)\subset \mathbb{R}$ for $i=2, 3,4,...,n$. Let $A$ be a quaternion algebra over $k$, which is ramified at all but the first real places, and thus there is an isomorphism
$$\sigma: A\otimes_\bbQ\mathbb{R}=\SL(2,\mathbb{R})\oplus \mathbb{H}\oplus \mathbb{H}\oplus\cdots\oplus\mathbb{H},$$
where $\mathbb{H}$ denotes Hamilton's quaternions. Denote $P$ the projection to the first factor.

Let $\mathcal{O}$ be an order in $A$, and $\mathcal{O}^1$  denote the group of elements of reduced norm 1. Then $P\sigma(\mathcal{O}^1)$ is a lattice of $\PSL(2,\mathbb{R})$. The class of
\emph{arithmetic Fuchsian groups} is all Fuchsian lattices commensurable with such groups $P\sigma(\mathcal{O}^1)$. In addition, we say that a Fuchsian group is \emph{derived from a
quaternion algebra} if it is a subgroup of finite index in some $P\sigma(\mathcal{O}^1)$.

\subsection{Characterization of arithmetic Fuchsian and Kleinian groups}
Takeuchi characterizes arithmetic Fuchsian groups within the class of all Fuchsian lattices in \cite{takeuchi1975characterization}. Maclachlan and Reid \cite{maclachlan1987commensurability} extend this work to Kleinian groups. 

Let $\Gamma$ be a Fuchsian group, and let $\Gamma^{(2)}$ denote the subgroup generated by the squares of elements of $\Gamma$. Note that if $\Gamma$ is finitely generated then $\Gamma^{(2)}$ is of finite index in $\Gamma$.
\begin{thm}[\cite{takeuchi1975characterization}, \cite{borel1981commensurability}]\label{arithmetic}
    If $\Gamma$ is an arithmetic Fuchsian group, then $\Gamma^{(2)}$ is derived from a quaternion algebra.
\end{thm}
\begin{thm}[\cite{takeuchi1975characterization}]\label{THM: Fuchsian}
Let $\Gamma$ be a cofinite Fuchsian group. Then $\Gamma$ is derived from a quaternion algebra over a totally real algebraic number field if and only if $\Gamma$ satisfies the following two conditions:
\begin{enumerate}
    \item 
$K:=\bbQ(\mathrm{Tr}(\Gamma))$ is an algebraic number field of finite degree and $\mathrm{Tr}(\Gamma)$
is contained in the ring of integers $\mathcal{O}_K$ of $K$.
\item
For any embedding $\phi$ of $K$ into $\mathbb{C}$, which is not the identity, $\phi(\mathrm{Tr}(\Gamma))$
is bounded in $\mathbb{C}$.
\end{enumerate}
\end{thm}

\subsection{$\mathbb{Z}$-GCD and $\mathbb{Z}$-LCM of real numbers}

We introduce the \textit{$\mathbb{Z}$-greatest common divisor} \textit{($\mathbb{Z}$-GCD)} and \textit{$\mathbb{Z}$-least common multiple} \textit{($\mathbb{Z}$-LCM)} in this subsection. In this paper, it is sufficient to consider only rational numbers. However, the theory extends naturally to all real numbers, and we will deal with the general setting. These two definitions are natural generalizations of their counterparts for integers to the $\mathbb{Z}$-module $\mathbb{R}$. It may also make sense for all $\mathbb{Z}$-modules, but this is beyond the scope of this paper.  

\begin{definition}
    Let $x$ and $y\in \mathbb{R}^+$. The \textit{$\mathbb{Z}$-GCD} of $x$ and $y$, denoted by $\GCD_{\mathbb{Z}}(x,y)$, is defined as follows:
    \begin{enumerate}
        \item If $\frac{x}{y}\notin \mathbb{Q}$, then $\GCD_{\mathbb{Z}}(x,y)=0$;
        \item If $\frac{x}{y}=\frac{p}{q}$ with $p$, $q\in \mathbb{Z}^+$ and $p$, $q$ coprime, then $\GCD_{\mathbb{Z}}(x,y)=\frac{x}{p}=\frac{y}{q}$.
    \end{enumerate}
\end{definition}

Similarly, we also have the following definition.

\begin{definition}
    Let $x$ and $y\in \mathbb{R}^+$. The \textit{$\mathbb{Z}$-LCM} of $x$ and $y$, denoted by $\LCM_{\mathbb{Z}}(x,y)$, is defined as follows:
    \begin{enumerate}
        \item If $\frac{x}{y}\notin \mathbb{Q}$, then $\LCM_{\mathbb{Z}}(x,y)=+\infty$;
        \item If $\frac{x}{y}=\frac{p}{q}$ with $p$, $q\in \mathbb{Z}^+$ and $p$, $q$ coprime, then $\LCM_{\mathbb{Z}}(x,y)=xq=py$.
    \end{enumerate}
\end{definition}

It is not hard to verify that some basic properties of the standard GCD and LCM still hold for these definitions. For example: if $0<\frac{x}{y}\in\mathbb{Q}$, then we have

\[\GCD_\mathbb{Z}(x,y)\cdot \LCM_\mathbb{Z}(xy)=xy.\]

\subsection{Natural density of $\mathbb{Z}$-affine space}\label{GCD}
Let $A\subset \mathbb{R}$. The \textit{natural density} of $A$ is defined by 
$$\rho(A)=\lim_{n\to \infty}\frac{\#\{A\cap [-n,n]\}}{2n},$$
whenever the limit exists.

A subset $A\subset \mathbb{R}$ is called a \textit{$\mathbb{Z}$-affine subspace} if there exist reals numbers $x$, $y>0$ such that $A=\{x+ky| k\in \mathbb{Z}\}.$ We denote this set as $A_{x,y}$. We also denote the set $\{x\}$ as $A_{x,\infty}$.  It is clear that $\rho(A_{x,y})=\frac{1}{y}$ for all $y$, where we take the convention $\frac{1}{\infty}=0$. 

We will apply the inclusion-exclusion argument later. Thus, we are also interested in the intersection of two $\mathbb{Z}$-affine subspaces. The intersection of $A_{x,y}$ and $A_{x',y'}$ has following possibilities:
\begin{enumerate}
    \item if $y=\infty$ or $y'=\infty$, then the intersection is 1 point or empty;
    \item if $y\neq \infty$ and $y'\neq \infty$, then:
    \begin{enumerate}
        \item if $\frac{y}{y'}\notin \mathbb{Q}$, then the intersection is a single point or empty;
        \item if $\frac{y}{y'}\in \mathbb{Q}$, then the intersection is either empty or a $\mathbb{Z}$-affine space $A_{x'',y''}$ with $y''=\LCM_{\mathbb{Z}}(y,y')$.
    \end{enumerate}
\end{enumerate}
Note that in all cases, we have $\rho(A_{x,y}\cap A_{x',y'})\leq \frac{1}{\LCM_{\mathbb{Z}}(y,y')}$.

\subsection{Dirichlet's Theorem on arithmetic progressions}
Dirichlet’s theorem on arithmetic progressions is a gem of number theory. A great part of its beauty lies in the simplicity of its statement.
\begin{theorem}[Dirichlet]
Let $a$, $m\in\mathbb{Z}$, with $\GCD(a, m)=1$. Then there are infinitely many prime numbers in the sequence of integers $a$, $a+m$, $a+2m$, $\cdots$, $a+km$, $\cdots$, for $k\in\mathbb{N}$.
\end{theorem}
We will need an effective version of Dirichlet's Theorem. For the same $a$ and $m$, denote the set
\[\Pi(x,m,a)=\{p\in \bbN \mid p\ \mathrm{prime},\  p\equiv a(\mathrm{mod}\ m),\ p\leq x\}.\]

In what follows, we consider two functions: $\pi(x,m,a)= \#\Pi(x,m,a)$ and
$$ S(x,m,a)=\sum_{p\in\Pi(x,m,a)}\frac{1}{p}-\frac{1}{\varphi(m)}\log\log x,$$
where $\varphi$ is Euler’s totient function. The following theorem is a slight variation of \cite{pomerance1977distribution}*{Theorem 1}. For completeness, we prove it here. The proof is very similar to the proof in \cite{pomerance1977distribution}.
\begin{theorem}\label{Dirichlet}
    There is an absolute constant $D$ such that for all $x\geq y \geq e^m$, all $a\geq 2$,  and $m$ with $\GCD(a,m)=1$, we have
    \[ |S(x,m,a)-S(y,m,a)|\leq \frac{D}{m\varphi(m)}.\]
\end{theorem}
\begin{proof}
    From the Siegel-Walfisz theorem (\cite{montgomery2007multiplicative}*{Corollary 11.19}), there is an absolute
constant $A$ such that for all $t>e^m$,
\begin{equation}\label{pi estimate}
    |\pi(t,m,a)-\frac{t}{\varphi(m)\log t}|< \frac{At}{\varphi(m)\log^2 t}
\end{equation}
An application of integration by parts yields
\begin{equation}\label{2.2}
    \begin{array}{rl}
    \sum_{p\in\Pi(x,m,a)\setminus \Pi(y,m,a)}\frac{1}{p}&=\int_y^x \frac{d\pi(t,m,a)}{t}\\&=\frac{1}{x}\pi(x,a,m)-\frac{1}{y}\pi(y,a,m)+\int_y^x \frac{\pi(t,m,a)}{t^2}dt.
    \end{array}
\end{equation}
Now by \eqref{pi estimate}, we have for all $t>e^m$,
\begin{equation}
    0 \leq \frac{1}{t}\pi(t,m,a)<\frac{1}{m\varphi(m)}+\frac{A}{m^2\varphi(m)}. 
\end{equation}
It follows that 
\begin{equation}\label{2.4}
    |\frac{1}{x}\pi(x,m,a)-\frac{1}{y}\pi(y,m,a)|<\frac{2}{m\varphi(m)}+\frac{2A}{m^2\varphi(m)}.
\end{equation}
Now since
\[
\int_y^x \frac{1}{\varphi(m)t\log t}dt=\frac{1}{\varphi(m)}(\log\log x-\log\log y).
\]
we have, by \eqref{pi estimate}
\begin{equation}\label{2.5}
   \begin{array}{rl} |\int_y^x \frac{\pi(t,m,a)}{t^2}dt-(\frac{1}{\varphi(m)}(\log\log x-\log\log y))|&<\int_y^x \frac{A}{\varphi(m)t\log^2t}dt \\
   &< \frac{A}{m\varphi(m)}.
   \end{array}
\end{equation}
Equations \eqref{2.2}, \eqref{2.4} and \eqref{2.5} imply that
\[
|S(x,m,a)-S(y,m,a)|<\frac{A+2}{m\varphi(m)}+\frac{2A}{m^2\varphi(m)}\leq \frac{3A+2}{m\varphi(m)}.
\]
\end{proof}

\subsection{Critical exponent and Cheeger constant}
Let $\Gamma$ be a Fuchsian group, and fix a point $o\in {\bfH}^2$. Consider the \textit{Poincar\'e series} of $\Gamma$, defined as:
\[ 
P_\Gamma(s)=\sum_{\gamma\in \Gamma}e^{-sd(o,\gamma o)},
\]
where $d$ is the hyperbolic metric on ${\bfH}^2$.
Then the \textit{critical exponent} of $\Gamma$ is given by
\[
\delta_\Gamma=\inf\{s\mid P_\Gamma(s)<\infty.\}
\]

For our purposes, the following aspect of the critical exponent suffices:
\begin{theorem}\cite{patterson1988lattice}\label{critical}
   Suppose $\Gamma$ is a convex cocompact subgroup of $\PSL(2,\mathbb{R})$. Then there exists $c>0$ and a $\Gamma$-invariant continuous function $F: {\bfH}^2\to \mathbb R^+$, such that 
   \[
   \#\{\gamma\in\Gamma\mid d(x,\gamma x')\leq R\}\sim cF(x)F(x')e^{\delta_\Gamma R}
   \]
   for all $x$, $x'\in{\bfH}^2$.
\end{theorem}

\medskip

Now we turn to the discussion of \textit{Cheeger constants}. In 1970, Cheeger \cite{cheeger1970lower} introduced an isoperimetric constant, now known as the Cheeger constant, to bound from below the first positive eigenvalue of the Laplacian. For any compact $n$-dimensional Riemannian manifold, the Cheeger constant of $M$ is given by
\[
h(M)=\inf\frac{\mathrm{Vol}(\partial A)}{\mathrm{Vol}(A)},
\]
where $A$ runs over all open subsets with $\mathrm{Vol}(A)\leq \frac{1}{2}\mathrm{Vol}(M)$. Cheeger \cite{cheeger1970lower} showed that
\[
\lambda_1(M)\geq \frac{1}{4}h^2(M),
\]
where $\lambda_1$ is the smallest positive eigenvalue of the Laplace-Beltrami operator. Cheeger's inequality also holds for non-compact Riemannian manifolds, provided that $A\cup\partial A$ is compact.

We also need a special case of a result by Buser.
\begin{theorem}\cite{buser1982note}*{Theorem 7.1}\label{buser}
    There exists a constant $\kappa$ such that for any non-compact hyperbolic surface $\Sigma$,
    \[
    \lambda_1(\Sigma)\leq \kappa h(\Sigma).
    \]
\end{theorem}
The final ingredient is the relationship between $\lambda_1$ of ${\bfH}^2/\Gamma$ and the critical exponent $\delta_\Gamma$ of a Fuchsian group $\Gamma$.

\begin{thm}\cites{elstrodt1973resolvente, elstrodt1973resolvente1, patterson1976limit}\label{cheeger}
    Let $\Gamma$ be a Fuchsian group, then 
    \[
    \lambda_1({\bfH}^2/\Gamma)=\left\{\begin{array}{lc}
    \frac{1}{4}&\delta_\Gamma<\frac{1}{2},\\
    \delta_\Gamma(1-\delta_\Gamma)& \delta_\Gamma\geq \frac{1}{2}.
    \end{array}
    \right.
    \]
\end{thm}
This result has been generalized to discrete subgroups of $\mathrm{PSO}(n,1)$ in \cite{corlette1990hausdorff}.
\section{second order linear recurrence sequences}
In this section, we consider a generalization of the Fibonacci sequence, which we call \textit{second order linear recurrence sequences (SLRS)}. All results in this section are crucial for Section~\ref{Step 2}. 

\begin{definition}
    A sequence of real numbers $F_n$, $n\geq 0$, is called a \textit{SLRS} if there exist reals $a$, $b$ such that $F_n$ satisfying the following recurrence relation for all $n\geq 2$:
    $$F_n=aF_{n-1}-bF_{n-2}.$$
\end{definition}

 We call $F_n$ a $SLRS(a,b)$. We will mainly consider the special case when all $F_n$ and $a$ are rational numbers and $b=1$. The following two lemmas will be used in Section~\ref{Step 2}.

Let $p$ and $q\geq 2$ be two positive integers and $\GCD(p,q)=1$. 

We have the following:
\begin{lemma}\label{estimate of boundedness}
Let $F_n$ be a non-constant $SLRS(\frac{p}{q},1)$. Denote the reduced form of $F_n$ by $\frac{f_n}{f_n'}$. Then there exist $j\geq 0$ and bounded sequence $i'_n\in \bbZ$ such that 
\[f_n'=i'_nq^{n-j}.\]
\end{lemma}
\begin{proof}
Let $M\in \mathbb{N}^+$ such that $MF_0$, $MF_1$ are integers. Define a new sequence $\overline{F_n}=Mq^nF_n$. Then $\overline{F_n}$ is a $SLRS(p,q^2)$. By the choice of $M$, $\overline{F_n}$ is an integer sequence. Define $H_n$ be an integer $SLRS(p,q^2)$ with $H_0\overline{F_1}-H_1\overline{F_0}\neq 0$. Note that this is possible since $\bar{F_0}$ and $\bar{F_1}$ can not both be 0; otherwise, $F_n$ will be a constant sequence.  

Define matrices $A_n=\begin{pmatrix}
    \overline{F_{n-1}}&H_{n-1}\\
    \overline{F_n}&H_n
\end{pmatrix}$
for all $n\in \mathbb N$. By the recurrence relation:
$$A_n=\begin{pmatrix}
    0& 1\\
    p& -q^2
\end{pmatrix}A_{n-1}=\begin{pmatrix}
    0& 1\\
    p& -q^2
\end{pmatrix}^nA_0.$$
Taking determinant, we have that $\GCD(\overline{F_{n-1}},\overline{F_n})$ is a factor of $p^n(\overline{F_0}H_1-\overline{F_1}H_0)$. 

Now for any prime factor $r$ of $q$, denote $\nu_r(n)=\max\{s\in \mathbb{N}\mid r^s|n\}$ as the evaluation at $r$. By the definition of $SLRS$, if $\nu_r(\overline{F_n})\leq \nu_r(\overline{F_{n-1}})$, then $\nu_r(\overline{F_{n+1}})=\nu_r(\overline{F_n})$. By induction,  $\nu_r(\overline{F_{n+k}})=\nu_r(\overline{F_n})$ for all $k\in \mathbb{N}$. Hence, the map $\mu_r: n\to v_r(\overline{F_{n}})$ is bounded unless it is a strictly increasing function. 

However, assume that $\mu_r$ is strictly increasing. Then $r^{n-1}\leq r^{v_r(\overline{F_{n-1}})}$ is a factor of $\GCD(\overline{F_{n-1}},\overline{F_n})$. The fact $\GCD(p,q)=1$ implies that $r^{n-1}$ is a factor of $\overline{F_0}H_1-\overline{F_1}H_0$ for all $n$, which is a contradiction. Hence, the map $\mu_r$ is bounded for all $r$. We conclude that $\GCD(\overline{F_n}, q^n)$ is bounded. 

Since $F_n=\frac{\overline{F_n}}{Mq^n}$, it follows that 
\[f_n=\frac{\overline{F_n}}{\GCD(\overline{F_n}, Mq^n)}\quad \mathrm{and}\quad f_n'=\frac{Mq^n}{\GCD(\overline{F_n}, Mq^n)}.\] 

Since
\[\GCD(\overline{F_n}, Mq^n) | M  \GCD(\overline{F_n}, q^n),\]
and there exists $j\geq 0$ such that $\GCD(\overline{F_n}, q^n)|q^j$, we have
$q^{n-j}|f_n'$, and $\frac{f_n'}{q^{n-j}}$ are bounded. Let $i'_n=\frac{f'_n}{q^{n-j}}$. The result follows.
\end{proof}
\begin{lemma}\label{estimate 2}
Let $F_n$ and $G_n$ be two $SLRS(\frac{p}{q},1)$, which are not the constant sequence $0$. Denote the reduced form of $F_n$, $G_n$ by $\frac{f_n}{f_n'}$, $\frac{g_n}{g_n'}$, respectively. 

    If $F_0G_1-F_1G_0\neq 0$, then $\GCD(f_n,g_n)$ is bounded.
\end{lemma}
\begin{proof}
 Let $M\in \mathbb{N}^+$ such that $MF_0$, $MF_1$, $MG_0$ and $MG_1$ are all integers. Define two new sequence $\overline{F_n}=Mq^nF_n$ and $\overline{G_n}=Mq^nG_n$. Then $\overline{F_n}$ and $\overline{G_n}$ are both $SLRS(p,q^2)$. By the choice of $M$, $\overline{F_n}$ and $\overline{G_n}$ are integer sequence. Let $H_n'=(\overline{G_1}-\overline{G_0})\overline{F_n}+(\overline{F_0}-\overline{F_1})\overline{G_n}$. The sequence $H_n'$ is also a $SLRS(p,q^2)$, and we have $H'_0=H'_1$.

For any prime factor $t$ of $p$, define $\nu_t$ as in the proof of Lemma~\ref{estimate of boundedness}. Since $H'_0=H'_1$, for any prime factor $t$ of $p$, $\nu_t(H_0)= \nu_t(H_1)$. Hence by $\GCD(t,q)=1$,
\[\nu_t(H_2)= \nu_t(pH_1-q^2H_0)=\nu_t(H_0).\]
 By induction, we conclude that $\nu_t(H'_n)=\nu_t(H'_0)$. Therefore, $\GCD(H'_n, p^n)$ is bounded. 

Since $H'_n$ is a linear combination of $\overline{F_n}$ and $\overline{G_n}$ with integer coefficients, it is clear that $\GCD(\overline{F_{n}},\overline{G_n})$ is a divisor of $H'_n$. 
 It follows that $\GCD(\overline{F_{n}},\overline{G_n}, p^n)$ is bounded.

On the other hand, define matrices \[A_n=\begin{pmatrix}
    \overline{F_{n}}&\overline{G_{n}}\\
    \overline{F_{n+1}}&\overline{G_{n+1}}
\end{pmatrix}\]
for all $n\in \mathbb N$. By the recurrence relation:
$$A_n=\begin{pmatrix}
    0& 1\\
    p& -q^2
\end{pmatrix}A_{n-1}=\begin{pmatrix}
    0& 1\\
    p& -q^2
\end{pmatrix}^{n}A_0.$$
Taking the determinant, we have that $\GCD(\overline{F_{n}},\overline{G_n})$ is a factor of $p^{n}(\overline{F_0}\overline{G_1}-\overline{F_1}\overline{G_0})$.  

The result follows from the estimation:
\[\GCD(\overline{F_{n}},\overline{G_n})\leq \GCD(\overline{F_{n}},\overline{G_n}, p^n)(\overline{F_0}\overline{G_1}-\overline{F_1}\overline{G_0})\] 
and the facts that:
\[f_n=\frac{\overline{F_n}}{\GCD(\overline{F_n}, Mq^n)}, \quad g_n=\frac{\overline{G_n}}{\GCD(\overline{G_n}, Mq^n)}.\]
Note that it has been shown that $\GCD(\overline{F_n}, Mq^n)$ and $\GCD(\overline{G_n}, Mq^n)$ are bounded in the proof of Lemma~\ref{estimate of boundedness}.
\end{proof}
\section{Proof of Theorem~\ref{main2}}
We prove Theorem~\ref{main2} in this section. 

Let $\Gamma<\SL(2,\bbR)$. Assume that $\Gamma$ contains a parabolic
element, and
\[
\#\{\mathrm{Tr}(\Gamma)\cap [-n,n]\}=o(n\log\log\log n).
\]

\subsection{Case 1: $\Gamma$ is elementary}\label{elementary}
Up to conjugacy, we may assume the $\Gamma$-fixed point on the boundary is $\infty$. Therefore $\Gamma$ is a subgroup of the Borel subgroup $B=\setdef{\begin{pmatrix}
    \tau &t\\
    0& \frac{1}{\tau}
\end{pmatrix}}{\tau\neq 0, t\in \mathbb{R}}.$ 

Let $\phi_{11}: \Gamma\to (\bbR^*,\cdot)$ be the group homeomorphism
\[
\phi_{11}(\begin{pmatrix}
    \tau &t\\
    0& \frac{1}{\tau}
\end{pmatrix})=\tau.
\]
Note that the image of $\phi_{11}$, denoted $D(\Gamma)$, is a subgroup of $(R^*,\cdot)$. Since $\mathrm{Tr}(\Gamma)$ is discrete, it follows that $D(\Gamma)$ is discrete as well. 

If $D(\Gamma)$ is finite, then $D(\Gamma)=\{\pm1\}\cong \bbZ/2\bbZ$ or trivial. In this case, set $\lambda=1$ and $\gamma=\begin{pmatrix}
    1& 0\\
    0& 1
\end{pmatrix}$.

If $D(\Gamma)$ is infinite, then $D(\Gamma)\cong \bbZ$ or $\bbZ\times \bbZ/2\bbZ$. Let $\lambda$ be a generator of the $\bbZ$-factor, and let $\gamma\in \phi_{11}^{-1}(\lambda)$. Since $\lambda\neq \pm 1$, up to conjugation, we may assume that $\gamma=\begin{pmatrix}
    \lambda& 0\\
    0&\frac{1}{\lambda}
\end{pmatrix}$. Note that the conjugation used here preserved the point $\infty$ on the boundary since it is the unique attracting fixed point of $\gamma$. Thus, it preserves the group $B$. In particular, the group $\Gamma$ is still upper-triangular.

The group $\Gamma$ acts on $\bbR^2\setminus\{(0,0)\}/\bbR^+$ by projective linear transformations. Let $\eta\in \Gamma$. The sign of $\phi_{11}(\eta)$, denoted $\epsilon(\eta)$, is determined by the relation $\eta([1:0])=[\epsilon(\eta)1: 0]$. Furthermore, $\phi^2_{11}(\eta)$ is simply the derivative $\eta'([1:0])$. Both quantities are invariants under conjugation by $g\in\SL(2,\bbR)$ that preserves $\infty$. Consequently, $\phi_{11}$is unchanged by the conjugation used to normalize $\gamma$.

Let $A=\langle\gamma\rangle$ and $Q=\phi_{11}^{-1}\{\pm 1\}$. Then it is clear $\Gamma=A\ltimes Q$. Since
\[
\begin{pmatrix}
    \lambda & 0\\
    0&\frac{1}{\lambda}
\end{pmatrix}\begin{pmatrix}
    \pm 1 & t\\
    0&\pm1
\end{pmatrix}\begin{pmatrix}
    \lambda & 0\\
    0&\frac{1}{\lambda}
\end{pmatrix}^{-1}=\begin{pmatrix}
    \pm 1 & \lambda^2 t\\
    0& \pm 1
\end{pmatrix},
\]
the set of upper-right entries of $Q$ is invariant under multiplication by $\lambda^2$.

Now let $P=\phi_{11}^{-1}\{1\}$, and $\phi_{12}:P\to (\bbR,+)$ be the group homomorphosm
\[
\phi_{12}\left(\begin{pmatrix}
    1& t\\
    0 & 1
\end{pmatrix}\right)=t.
\] Then the image of $\phi_{12}$, $R$, is a subgroup of $(\bbR,+)$. 

This completes the proof of Theorem~\ref{main2} in the elementary case.
\subsection{Case 2: $\Gamma$ is non-elementary} The next few subsections will concentrate on the proof of the case of a non-elementary subgroup.

The first part of the proof is very similar to the work in \cite{hao2023bounded}. We provide more details here for better readability.

First, let us find a canonical form of an embedding that will significantly reduce the computations. 

 Let $\Gamma$ be a subgroup of $\SL(2, \mathbb R)$ with parabolic elements. Up to finite index, we identify $\Gamma$ with its image in $\PSL(2,\bbR)$. Let $x\in \partial\Gamma$ be a fixed point of a parabolic element with a parabolic subgroup $\Gamma_x$. Taking $g\in \Gamma$ with $g\cdot x\neq x$, then $g\cdot x$ is a parabolic fixed point of $\Gamma$, and the parabolic subgroup is given by $\Gamma_{g\cdot x}=g\Gamma_x g^{-1}$. Up to conjugation in $\PSL(2,\mathbb{R})$, we may assume $x=\infty$, $g\cdot x=0$. Since $g\cdot \infty=0$, $g$ is in the form $\begin{pmatrix}
    0&\frac{1}{\beta_0}\\
    -\beta_0&*
\end{pmatrix}$ for some $\beta_0\in\bbR$. If \[\Gamma_\infty=\setdef{\begin{pmatrix}
    1&k\\
    0&1
\end{pmatrix}}{k\in P}\] for a non-trivial subgroup $P$ of $(\mathbb{R},+)$, then a direct computation shows that \[\Gamma_0=\setdef{\begin{pmatrix}
    1&0\\
    k\beta_0^2&1
\end{pmatrix}}{k\in P}.\]

If $P$ is dense, then since $\begin{pmatrix}
    1&k\\
    0&1
\end{pmatrix}\times \begin{pmatrix}
    1&0\\
    l\beta_0^2&1
    \end{pmatrix}=\begin{pmatrix}
    1+kl\beta_0^2&k\\
    l\beta_0^2&1
    \end{pmatrix},$ the trace set $\mathrm{Tr}(\Gamma)$ will contains a dense subset of $\bbR$, a contradiction. It follows that $P$ is discrete. Note that for all $\lambda\in D(\Gamma)$, where $D(\Gamma)$ is defined in section~\ref{elementary}, $P$ is $\lambda^2$-invariant. Hence, $D(\Gamma)$ is finite. And $D(\Gamma)$ is trivial after we identify $\Gamma$ with its image in $\PSL(2,\bbR)$.
    
    By taking conjugation by $\begin{pmatrix}
    \lambda&0\\
    0&\frac{1}{\lambda}
    \end{pmatrix}$ for suitable $\lambda$, we may take $P=\bbZ$ and denote $\beta=\frac{\beta_0^2}{\lambda^2}$, which is the left corner of $g$ after this conjugation. Hence we may assume $$\Gamma_\infty=\setdef{\begin{pmatrix}
    1&k\\
    0&1
\end{pmatrix}}{k\in \bbZ}.$$

Since $g\cdot \infty=0$, $g$ is in the form $\begin{pmatrix}
    0&\frac{1}{\beta}\\
    -\beta&*
\end{pmatrix}$. It follows from the fact $\Gamma_0=g\Gamma_\infty g^{-1}$ that
$$\Gamma_0=\setdef{\begin{pmatrix}
    1&0\\
    k\beta^2&1
\end{pmatrix}}{k\in \bbZ}.$$

From now on to the end of this section, we assume $\Gamma$ is a subgroup of $\SL(2,\bbR)$ containing $\Gamma_0$ and $\Gamma_\infty$ as above in this section. We will show that, under this embedding, $\Gamma$ is commensurable with a subgroup of $\PSL(2,\bbZ)$ when the growth rate of the trace set is slower than $O(n\log\log\log n)$. 

To prove it, we first show that $\Gamma^{(2)}<\PSL(2,\bbQ)$ when the growth of the trace set less than $n\log n$ in Section~\ref{A^2 section}. Then we proceed by contradiction in Section~\ref{Step 2}, that for any element whose trace is not in $\mathbb{Z}$, we construct a family of elements with trace set grows at least $n\log\log\log n$. In Section~\ref{Finite}, we show that the quotient group $\Gamma/(\Gamma\cap \SL(2,\mathbb Q))$ is finite. And in Section~\ref{Commen}, by a similar strategy to that of Takeuchi, we show that $\bar{\Gamma}:=\Gamma\cap \SL(2,\mathbb Q)$ is commensurable to a subgroup of $\SL(2,\bbZ)$.

\subsection{First step: $\Gamma^{(2)}$ is rational.}\label{A^2 section}
The main result in this subsection is Lemma~\ref{A^2}. 

 Let $x_1,x_2,\cdots,x_n\in\bbR$. Denote $\bbQ\langle x_1,x_2,\cdots,x_n\rangle$ the $\bbQ$-vector space generated by $x_1$, $x_2,\cdots, x_n$. 
\begin{lemma}\label{A^2}
If the growth of $\mathrm{Tr}(\Gamma)$ is $o(n\log n)$, then $A^2\in \PSL(2,\bbQ)$ for all $A\in \Gamma$.
\end{lemma}
\begin{proof}
Let $A=\begin{pmatrix}
    a&b\\c&d
\end{pmatrix}\in \Gamma$.
Without loss of generality, we assume $c\neq 0$. Indeed, $c=0$ implies that $A\cdot \infty=\infty$. Then $A\in \Gamma_\infty\subset \PSL(2,\bbQ)$.

The key step to prove Lemma~\ref{A^2} is the following;
\begin{claim}\label{rank 1}
    $\bbQ\langle\beta^2a,\beta^2b,c,\beta^2d\rangle=\bbQ\langle c\rangle.$
\end{claim}

The proof of Claim~\ref{rank 1} has two steps. Both are based on the analysis of a particular subset of the trace set. We construct the set first.

Since
$$\begin{pmatrix}
  1&0\\
  k\beta^2&1
\end{pmatrix}\begin{pmatrix}
    a&b\\
    c&d
\end{pmatrix}\begin{pmatrix}
    1&l\\
    0&1
\end{pmatrix}=\begin{pmatrix}
    a&al+b\\
    k\beta^2a+c&kl\beta^2a+lc+k\beta^2b+d
\end{pmatrix},$$
$\mathrm{Tr}(\Gamma)$ contains all elements of the form $a+d+kl\beta^2a+k\beta^2b+lc$, $(k,l)\in \bbZ^2$. 

Hence the set
\[
    \Omega_{a,b,c}:=\setdef{kl\beta^2a+k\beta^2b+lc}{(k,l)\in\bbZ^2}
\]
has growth less than $O(n\log n)$. For convenience, denote $\Theta(k,l)=kl\beta^2a+k\beta^2b+lc$.

\textbf{Step 1.} $\beta^2 a\in \bbQ\langle c,\beta^2b\rangle.$ 

We prove this by contradiction. Assume $\beta^2 a\notin \bbQ\langle c,\beta^2b\rangle$. Let $(k,l)\in \mathbb{Z}^2$ with $kl\neq 0$. Consider the equation for $(k',l')\in \mathbb{Z}^2$: $\Theta(k,l)=\Theta(k',l')$. Then $kl=k'l'$ and $k\beta^2 b+lc=k'\beta^2 b+l'c.$ Hence $l'=\frac{kl}{k'}$, and we have
$$(k-k')\beta^2 b=(\frac{l(k-k')}{k'})c.$$

\begin{enumerate}
    \item If $k'=k$, then $l'=l$. 
    \item If $k'\neq k$, then $\frac{l}{k'}=\frac{\beta^2 b}{c}$. And therefore $k'=\frac{lc}{\beta^2 b}$, $l'=\frac{k\beta^2 b}{c}$.
\end{enumerate}
Therefore $\Theta$ is at most 2 to 1 on the set $\{(k,l)\in \mathbb{Z}^2\mid kl\neq 0\}$.

Considering the set 
$$D_N:=\setdef{(k,l)}{1\leq kl\leq N,k,l\in \bbN^+},$$ 
it has $\sum_{j=1}^N \lfloor\frac{N}{j}\rfloor\geq N\ln N-N$ many elements. And all elements in $\Theta(D_N)$ have absolute value less than $N(|\beta^2 a|+|\beta^2 b|+|c|)$. Therefore, the trace set grows at least by $O(n\log n)$, a contradiction.  
Therefore $\beta^2 a\in \bbQ\langle c,\beta^2b\rangle.$

\textbf{Step 2.} $\beta^2 b\in \bbQ\langle c\rangle.$ 

Since $\beta^2 a\in \bbQ\langle c,\beta^2b\rangle,$ there exist $s,t\in\bbQ$ with $\beta^2 a=s\beta^2 b+t c$. Now the set
\[\Omega_{a,b,c}=\setdef{(skl+k)\beta^2b+(tkl+l)c}{k,l\in\bbZ}\]
has growth $o(n\log n)$. 

 Define $$\Phi(k,l)=(skl+k,tkl+l).$$ And notice that $$\Theta(k,l)=(skl+k)\beta^2b+(tkl+l)c.$$

 We consider all possibilities of the pair $s$, $t$. 
\begin{enumerate}
    \item Case 1: $s=t=0$. Then $\Theta$ map the set $\setdef{(k,l)}{1\leq l,k\leq N}$ to numbers with norm no more that $N(|\beta^2 b|+|c|)$. And $\Phi$ is injective.
    \item Case 2: $s=0$, $t\neq 0$.
 The image of $\Phi$ determines $k$. It follows that $\Phi$ is injective when $k\neq -\frac{1}{t}$. The set $D_N\setminus \{k=-\frac{1}{t}\}$ has more than $N\ln N-2N$ elements. And $|\Theta(u)|\leq N[(|s|+1)|\beta^2 b|+(|t|+1)|c|]$ for all $u\in D_N\setminus \{k=-\frac{1}{t}\}$. 
     \item Case 3: $s\neq 0$ and $t=0$. Similar to Case 2.
     \item Case 4: $st\neq 0$. Then $|\Theta(u)|\leq N[(|s|+1)|\beta^2 b|+(|t|+1)|c|]$ for all $u\in D_N$. And $\Phi$ is at most 2 to 1.

     To see that $\Phi$ is at most 2 to 1 in this case, fix $(k,l)\in D_N$. Assume that $\Phi(k,l)=\Phi(k',l')$. Then we have $tk-sl=tk'-sl'$, hence $l'=\frac{t}{s}(k'-k)+l.$ Plugging this into $skl+k=sk'l'+k'$, we have a quadratic equation for $k'$. It is clear $k'=k$ is a solution, hence the other solution is $k'=\frac{-(sl+1)}{t}$. In conclusion, there are at most two solutions in $D_N$ with $\Phi(K,l)=\Phi(k',l')$.
\end{enumerate}

In all cases, $\Theta$ maps a set with growth at least $O(n\log n)$ to a set of growth $n$ (up to a constant), and $\Phi$ is finite to one on this set. Dirichlet's Pigeonhole principle gives  $(k,l)\neq(k',l')$ such that $\Theta(k,l)=\Theta(k',l')$ and $\Phi(k,l)\neq\Phi(k',l')$ for $N$ big enough. We have a nontrivial homogeneous linear equation of $\beta^2b$ and $c$. If $b=0$, step 2 is trivially true. If $b\neq 0$, and by assumption, $c\neq 0$, we have $\beta^2 b\in \bbQ\langle c\rangle.$ 

Similarly, $\beta^2d\in \bbQ\langle c\rangle$ by a similar consideration on the following family of elements:
\[
\begin{pmatrix}
    1&l\\
    0&1
\end{pmatrix}\begin{pmatrix}
    a&b\\
    c&d
\end{pmatrix}\begin{pmatrix}
    1&0\\
    k\beta^2&1
\end{pmatrix}=\begin{pmatrix}
    a+kl\beta^2d+lc+k\beta^2b&ld+b\\
    k\beta^2d+c&d
\end{pmatrix}.
\]

The claim is proved.

Now we prove the lemma.

Considering the element $\begin{pmatrix}
    1&1\\
    0&1
\end{pmatrix}\begin{pmatrix}
    1&0\\
    \beta^2&1
\end{pmatrix}=\begin{pmatrix}
    1+\beta^2&1\\
    \beta^2&1
\end{pmatrix}\in\Gamma$.
Claim~\ref{rank 1} gives $\beta^2+\beta^4\in \bbQ\langle \beta^2\rangle$. We conclude that $\beta^2\in \bbQ$. Then 
\[\bbQ\langle a,b,c,d\rangle=\bbQ\langle c\rangle.\]

Now $A=cA'$ with $A'\in \GL(2,\bbQ)$. Taking determinant, $c^2\in \bbQ$. Finally,
\[A^2=c^2A'^2\in \PSL(2,\bbQ).\] 
\end{proof}

\subsection{Second step: $\mathrm{Tr}(\bar{\Gamma})$ is in $\bbZ$.}\label{Step 2}

By Theorem~\ref{arithmetic} and Theorem~\ref{THM: Fuchsian}, it is enough to work with $\Gamma^{(2)}$ to show that $\Gamma$ is arithmetic when $\Gamma$ is a lattice. However, in general, we will work with a slightly larger subgroup $\bar{\Gamma}=\Gamma\cap \PSL(2,\bbQ)$. Now $\Gamma^{(2)}$ is a normal subgroup. And $\Gamma/\Gamma^{(2)}$ is a 2-group, hence a commutative group. We know that $\bar{\Gamma}$ is a normal subgroup of $\Gamma$. 

 We begin with a few simple lemmas.
 \begin{lemma}\label{Z-affine}
 If $A=\begin{pmatrix}
     a&b\\
     c&d
 \end{pmatrix}\in \bar{\Gamma}$, then $\mathrm{Tr}(\bar{\Gamma})$ contains the following subset:
 \[\setdef{a+d+k\beta^2b}{k\in\mathbb{Z}}\].
 \end{lemma}
 \begin{proof}
     This follows since \[A\begin{pmatrix}
         1&o\\
         k\beta^2&1
     \end{pmatrix}=\begin{pmatrix}
         a+k\beta^2 b&b\\
         c+k\beta^2 d&d
     \end{pmatrix}\in\bar{\Gamma}.\]
 \end{proof}
\begin{lemma}\label{left coner}
     If $A=\begin{pmatrix}
     a&b\\
     c&d
 \end{pmatrix}\in \bar{\Gamma}$, then 
 \[A\begin{pmatrix}
     1&l\\
     0&1
 \end{pmatrix}=\begin{pmatrix}
     a& al+b\\
     c& cl+d
 \end{pmatrix}\in\bar\Gamma\] for all $l\in \mathbb Z$.
\end{lemma}

Now we are ready to prove the key result of this section.
\begin{lemma}\label{Z}
 If $\mathrm{Tr}(\bar{\Gamma})$ has growth $o(n\log\log\log n)$, then $\mathrm{Tr}(\bar{\Gamma})$ is a subset of $\mathbb{Z}$.
\end{lemma}
\begin{proof}
    We prove this by contradiction. Assume that there exists $A=\begin{pmatrix}
    a&b\\
    c&d
\end{pmatrix}\in \bar{\Gamma}$ with trace $a+d\notin \mathbb{Z}$. Denote $a+b=\frac{p}{q}$ with $p,q\in\mathbb N$ and $p$, $q$ coprime to each other. 

First, note that we have $bc\neq 0$, since $A\notin \bar{\Gamma}_{\infty}$ or $\bar{\Gamma}_{0}$. By Lemma~\ref{Z-affine}, and up to multiplying by $-1$, we may assume trace $a+d>2$ and $d<0$. 

Let $A^n=\begin{pmatrix}
    a_n&b_n\\
    c_n&d_n, 
\end{pmatrix}$. By Cayley–Hamilton theorem, $A^2-\frac{p}{q}A+I=0$. Hence $A^n-\frac{p}{q}A^{n-1}+A^{n-2}=0$. 
It follows that the four sequences $a_n$, $b_n$, $c_n$ and $d_n$ are $SLRS(\frac{p}{q}, 1)$. Let $\lambda>1$ and $\frac{1}{\lambda}$ be the solution of the quadratic equation $x^2-\frac{p}{q}x+1=0$. Then there exist $\alpha$, $\beta$ so that 
\begin{equation}\label{SLRS}
a_n=\alpha \lambda^n+\frac{\beta}{\lambda^n}
\end{equation}
by the general theory of (homogeneous) linear recurrence sequences \cite{andrica2020recurrent}*{Section 2.2.1}. Set $n=0$ in equation \eqref{SLRS}, and since $a_0=1$, we have $\alpha+\beta=1$. Set  $n=1$ in equation \eqref{SLRS}, and note that by assumption $a_1=a>2$, we have $\alpha>1$. Therefore, $a_n>0$ and is increasing.

By Lemma~\ref{Z-affine} and \ref{left coner},  $\mathrm{Tr}(\bar{\Gamma})$ contains the family of $Z$-affine subspaces 
\[A(n,l):=A_{a_n+l c_n+d_n,\beta^2(l a_n+b_n)}\] for all $n$,$l\in \mathbb{Z}$. Note that $$\rho(A(n,l))=\frac{1}{\beta^2(l a_n+b_n)},$$ $$\rho(A(n,l)\cap A(n',l'))\leq \frac{1}{\beta^2\LCM_\bbZ(l a_n+b_n,l' a_{n'}+b_{n'})}.$$

We will select a subfamily such that the density of their union is infinite. For this, denote the reduced rational representation of $a_n$, $b_n$ by $\frac{A_n}{A'_n}$ and $\frac{B_n}{B'_n}$, respectively. Then \[l a_n+b_n=\frac{K_n(l S_n+T_n)}{L_n}\] where, $L_n=\LCM(A'_n,B'_n)$, $\GCD(S_n,T_n)=1$, $K_n=\GCD(\frac{A_nB_n'}{\GCD(A'_n,B'_n)}, \frac{B_nA_n'}{\GCD(A'_n,B'_n)})$. By Lemma~\ref{estimate of boundedness}, the reduced form of $\frac{L_n}{q^n}$ has a bounded numerator and denominator. In particular, 
\begin{equation}\label{L_n}
L_n\sim q^n.
\end{equation}  Also by Lemma~\ref{estimate of boundedness}, $\frac{B'_n}{\GCD(A'_n,B'_n)}$ and $\frac{A_n'}{\GCD(A'_n,B'_n)}$ are bounded. By Lemma~\ref{estimate 2} and 
\[
K_n\leq \GCD(A_n, B_n) \frac{B'_n}{\GCD(A'_n,B'_n)}\frac{A_n'}{\GCD(A'_n,B'_n)},
\] $K_n$ is bounded. On one hand, $S_n=L_n\frac{A_n}{A_n'K_n}\sim q^n a_n$. On the other hand, by equation~\eqref{SLRS}, $a_n\sim \lambda^n$. We conclude that 
\begin{equation}\label{S_n}
S_n\sim (\lambda q)^n.
\end{equation}

Now we further normalize the sequence, up to changing the matrix $A$ to $A^k$, we may assume that $S_n\geq 3$ for all $n\in \bbN^+$.

Define a new sequence $E_n=\exp(\exp(e_n))$, where: 
\[e_n=\sum_{i=0}^n \frac{\beta^2K_i\varphi(S_i)}{L_i}.\] 
Let 
\[I_n=\setdef{l\in\mathbb{Z}}{lS_n+T_n \mathrm{\ is\  prime}, E_{n-1}< l S_n+T_n\leq E_n}.\]
Finally, the family of $\mathbb{Z}$-affine subspaces is given by all $\{A(n,l)|n\geq 2,\ l\in I_n\}$. Let $U$ be the union of all such $A(n,l)$.

We continue with some estimations to finish the proof.

 \textbf{(I).} It is known (for example, see \cite{rosser1962approximate}) that for $n>2$,
$$\varphi(n)>\frac{n}{e^{\gamma}\log\log n+\frac{3}{\log\log n}},$$
where $\gamma$ is the Euler constant. It leads to
\[\frac{S_n}{L_n(e^\gamma\log\log S_n+\frac{3}{\log\log S_n})}\leq \frac{\varphi(S_n)}{L_n} \leq \frac{S_n}{L_n}.\] 

By equation~\eqref{L_n} and \eqref{S_n}, for $n\geq 2$,
\[\frac{S_n}{L_n(e^\gamma\log\log S_n+\frac{3}{\log\log S_n})}\sim \frac{\lambda^n}{\log n},\]
and 
\[
\frac{S_n}{L_n}\sim \lambda^n
\]
Thus, for $n\geq 2$, 
\[ \frac{\lambda^n}{\log n}=O(\frac{\varphi(S_n)}{L_n})\  \mathrm{and}\  \frac{\varphi(S_n)}{L_n}=O(\lambda^n).\]

Therefore, there exist $A>0$ and $B>0$ such that for all $n>2$,  
\[
e_n=\sum_{i=0}^n \frac{\beta^2K_i\varphi(S_i)}{L_i}\geq A\frac{\lambda^n}{\log n},
\]
and 
\[
e_n=\sum_{i=0}^n \frac{\beta^2K_i\varphi(S_i)}{L_i}\leq \sum_{i=0}^n B\lambda^n=O(\lambda^n).
\]

Hence $\log\log\log E_n=\log(e_n)\sim n$.

\textbf{(II).} We estimate $\sum_{l\in I_n, n\geq 2} \rho(A(n,l))$ here.

First, by definition of natural density in Section 2.D,
\[
\sum_{l\in I_n} \rho(A(n,l))=\frac{L_n}{\beta^2 K_n}\sum_{l\in I_n} \frac{1}{lS_n+T_n}
\]
Since $\log\log\log(E_n)\sim n$, for $n$ sufficiently large, $E_{n-1}\geq e^{S_n}$.
Applying Theorem~\ref{Dirichlet} with $x=E_n$, $y=E_{n-1}$, $m=S_n$ and $a=T_n$, for $n$ sufficiently large, we obtain
\[
|\sum_{l\in I_n} \frac{1}{lS_n+T_n}-\frac{1}{\varphi(S_n)}(e_n-e_{n-1})|\leq \frac{D}{S_n\varphi(S_n)}.
\]
Thus
\[
|\sum_{l\in I_n} \rho(A(n,l))-1|\leq \frac{DL_n}{\beta^2K_nS_n\varphi(S_n)}.
\]
By equation~\eqref{L_n}, \eqref{S_n}, and the fact $K_n$ is bounded, 
\[
\lim_{n\to \infty} \frac{DL_n}{\beta^2K_nS_n\varphi(S_n)}=0.
\]
Therefore, for $n$ sufficiently large, 
\[
\frac{1}{2}\leq \sum_{l\in I_n} \rho(A(n,l))\leq \frac{3}{2}.
\]
Thus,
\[\sum_{l\in I_n, n\geq 2} \rho(A(n,l))=\infty.\]

Also Let
\[
C^+=\sup_n\{\sum_{l\in I_n} \rho(A(n,l))\}>0,\]
\[
C_-=\inf_n\{\sum_{l\in I_n} \rho(A(n,l))\}>0.
\]

\textbf{(III).} Now for $(n,l)\neq (n',l')$, $$\rho(A(n,l)\cap A(n',l'))\leq \frac{1}{\beta^2\LCM_\bbZ(l a_n+b_n,l' a_{n'}+b_{n'})}.$$
Since \[\LCM_\bbZ(l a_n+b_n,l' a_{n'}+b_{n'})=\frac{\LCM(K_n, K_{n'})(lS_n+T_n)(l'S_{n'}+T_{n'})}{\GCD(L_n, L_{n'})}\] for sufficiently large $n$, $n'$, (so that $lS_n+T_n$, $l'S_{n'}+T_{n'}$ are larger compare to all of the four numbers: $K_n$, $K_{n'}$, $L_n$ and $L_{n'}$. This is possible since $q^n=o(E(n))$), $\LCM_\bbZ(l a_n+b_n,l' a_{n'}+b_{n'})$ is equal to $(l a_n+b_n)(l' a_{n'}+b_{n'})q^{\max\{n,n'\}}$ up to a constant factor. Hence, for $n\neq n'$,
\[\sum_{l\in I_n,\  l'\in I_{n'}}\rho(A(n,l)\cap A(n',l'))\leq C_1\sum_{l\in I_n,\  l'\in I_{n'}}\frac{\rho(A(n,l))\rho(A(n',l'))}{q^{\max\{n,n'\}}} \leq \frac{C_1(C^+)^2}{q^{\max\{n,n'\}}},\]
and for $n=n'$,
\[\sum_{l\neq l'\in I_n}\rho(A(n,l)\cap A(n,l'))\leq C_1\sum_{l\neq l'\in I_n}\frac{\rho(A(n,l))\rho(A(n,l'))}{q^{n}} \leq \frac{C_1(C^+)^2}{q^{n}}.\]
Therefore,
\[\sum_{l\in I_n,\  l'\in I_{n'},\  (n,l)\neq (n', l')}\rho(A(n,l)\cap A(n', l'))\leq \sum_n C_1(C^+)^2\frac{n}{q^n}\leq C_2.\]

\textbf{(IV).} Fix a number $k>e^{10}$, and let $\vartheta(k)=\max\{n\mid E_n\leq k^{1/4}\}$. By estimation (I),  $\vartheta(k)\sim \log\log\log k$. 

We will calculate the contribution of $U\cap [0,k]$ from all $A(n,l)$ with $n\leq \vartheta(k)$. For each such $A(n,l)$, for large enough $k$, we have: 
\[\#\{A(n,l)\cap [0,k]\}\geq \frac{k}{\beta^2(la_n+b_n)}-1\geq \frac{k\rho(A(n,l))}{2}\]
since $\frac{k}{\beta^2(la_n+b_n)}\geq O(\frac{kq^{n}}{k^{1/4}})>2$.

Thus:
\[ \sum_{l\in I_n, \ n\leq \vartheta(k)} \#\{A(n,l)\cap [0,k]\}\geq \frac{C_-k\vartheta(k)}{2}\sim \frac{C_- k\log\log\log k}{2}.\]

On the other hand, for large $k$, for two different set $A(n,l)$ and $A(n',l')$ with $n$, $n'\leq \vartheta(k)$, we have:
\[
\begin{array}{rl}
     \#\{A(n,l)\cap A(n',l')\cap [0,k]\}&\leq \frac{k}{\beta^2\LCM_\bbZ(la_n+b_n,l'a_{n'}+b_{n'})}+1  \\
     & \leq 2k\rho(A(n,l)\cap A(n',l'))
\end{array}\]
since $\frac{k}{\beta^2\LCM_\bbZ(la_n+b_n,l'a_{n'}+b_{n'})}\geq O(\frac{kq^{\max\{n,n'\}}}{k^{1/2}})>2$.

Therefore, 
\[ \sum_{l\in I_n,\ l'\in I_{n'} \ n, n'<\vartheta(k),\ (n,l)\neq (n',l')} \#\{A(n,l)\cap A(n',l')\cap [0,k]\}\leq 2C_2k.\]

By the Bonferroni inequality, also known as the inclusion-exclusion inequality, for example, taking the form in \cite{dohmen1999improved}*{Introduction} with $n=2$, and integration with respect to the counting measure on $\{U\cap [0,k]\}$, we get that there exists $C>0$ such that
\[ \#\{U\cap [0,k]\}\geq k[\frac{C\log\log\log k}{2}-2C_2].\]
This is a contradiction. Therefore, $\mathrm{Tr}(\bar{\Gamma})$ is a subset of $Z$.
\end{proof}

\subsection{Third step: $\Gamma/\bar{\Gamma}$ is finite.}\label{Finite}
The fact $\mathrm{Tr}(\bar{\Gamma})\subset \mathbb Z$ has strong consequences for the structure of $\bar{\Gamma}$. We will show this first. The results of this and the next subsections are deduced from the structure of $\bar{\Gamma}.$

\begin{lemma}\label{Structure lemma}
    There exists $N\in \bbN^+$ such that
    \[
    \bar{\Gamma}\subset \setdef{\begin{pmatrix}
        a&b\\
        c&d
    \end{pmatrix}}{Na\in \bbZ, Nb\in\bbZ,c\in\bbZ,Nd\in \bbZ}.
    \]
\end{lemma}
\begin{proof}
    Let $A=\begin{pmatrix}
        a&b\\
        c&d
    \end{pmatrix}\in\bar{\Gamma}$. Then $a+d\in \bbZ.$
    
    \textbf{(I).} Considering $A\begin{pmatrix}
        1&1\\
        0&1
    \end{pmatrix}=\begin{pmatrix}
        a&a+b\\
        c&c+d
    \end{pmatrix}$. We have $a+d+c\in\bbZ$, so $c\in Z$. In particular, $\beta^2\in\bbZ$. 
    Let $N=\beta^2$.

    \textbf{(II).} Considering $A\begin{pmatrix}
        1&0\\
        N&1
    \end{pmatrix}=\begin{pmatrix}
        a+Nb&b\\
        c+Nd&d
    \end{pmatrix}$. We have $a+Nb+d\in \bbZ$ and $c+Nd\in Z$. It follows that $Nb\in \bbZ$ and $Nd\in \bbZ$.

    \textbf{(III).} Finally, $a+d\in\bbZ$. Thus, $Na\in\bbZ$.

    This completes the proof.
\end{proof}
Now we are ready to show the main result of this subsection.
\begin{lemma}
    The quotient group $\Gamma/\bar{\Gamma}$ is finite.
\end{lemma}
\begin{proof}
     Let $A=\begin{pmatrix}
        a&b\\
        c&d
    \end{pmatrix}\in \Gamma$ with $c\neq 0$.  Then $A=cB$ with $B\in\SL(2,\mathbb Q)$, by the proof of Lemma~\ref{A^2}. Since $c^2\in\mathbb{Q}$, there is a unique square-free $D_A\in \bbN^+$ such that $A=\sqrt{D_A}B'$ with $B'\in\GL(2,\mathbb Q)$. For elements with $c=0$, we take $D_A=1$. 
    
    Then the map 
    \[D_\Gamma:\Gamma/\bar\Gamma\to \bbN^+/(\bbN^+)^2=\bigoplus_{p\ \mathrm{prime}} \bbZ/2\bbZ\]
    given by $D_\Gamma(A\bar\Gamma)=D_A$, is a well-defined and injective group homomorphism. Hence, it is sufficient to show that the image of $D_\Gamma$ is finite. 
    
    Let $A\in \Gamma$ with $A=\sqrt{D_A}\begin{pmatrix}
       b_1&b_2\\
       b_3&b_4
    \end{pmatrix}$. Then $b_1b_4-b_2b_3=\frac{1}{D_A}$, and the inverse of $A$ is $A^{-1}=\sqrt{D_A}\begin{pmatrix}
       b_4&-b_2\\
       -b_3&b_1
       \end{pmatrix}$. Since $\Bar\Gamma$ is a normal subgroup of $\Gamma$, we have:
       \[A\begin{pmatrix}
           1&1\\
           0&1
       \end{pmatrix}A^{-1}=\begin{pmatrix}
           1-Db_1b_3&Db_1^2\\
           -Db_3^2&1+Db_1b_3
           \end{pmatrix}\in\bar{\Gamma},\]
           and 
           \[ A\begin{pmatrix}
           1&0\\
           N&1
       \end{pmatrix}A^{-1}=\begin{pmatrix}
           1-NDb_2b_4&-DNb_2^2\\
           -DNb_4^2&1+DNb_2b_4
           \end{pmatrix}\in\bar{\Gamma}.\]

    By Lemma~\ref{Structure lemma}, and the fact that $D$ is square-free, the following is true:
    \begin{enumerate}
        \item $Db_3^2\in \bbZ$, hence $b_3\in\bbZ$,
        \item $DNb_1^2\in\bbZ$, hence $Nb_1\in\bbZ$,
        \item $DN^2 b_4^2\in\bbZ$, hence $Nb_4\in\bbZ$,
        \item $DN^2 b_2^2\in\bbZ$, hence $N b_2\in \bbZ$.
    \end{enumerate}
    In particular, the denominator of $\frac{1}{D_A}=b_1b_4-b_2b_3$ is a factor of $N^2$. Hence, $D_A$ has only finitely many choices, and it follows that the image of $D_\Gamma$ is finite. 
    
    The proof is completed.
\end{proof}

\subsection{Fourth step: $\bar{\Gamma}$ has a finite subgroup in $\SL(2,\mathbb{Z})$}\label{Commen} 
The proof here is essentially the same as in Takeuchi's work \cite{kisao1969some}.  However, in our case, the argument can be made more elementary.

\begin{lemma}
    The group $\bar{\Gamma}$ has a finite index subgroup which is also a subgroup of $\SL(2,\mathbb{Z}).$
\end{lemma}
\begin{proof}
    Consider the group ring $\bbZ[\bar\Gamma]$. Since  $\bar{\Gamma}$ is a subgroup of $\SL(2,\bbQ)$, there is a natural map $\varrho: \bbZ[\bar\Gamma]\to M(2,\bbQ)$ where $M(2,\bbQ)$ is the set of $2\times 2$ matrices over $\bbQ$. 
    
    It is sufficient to show that $\varrho( \bbZ[\bar\Gamma])$ is an order in $M(2,\bbQ)$, since $\bar\Gamma$ is a subgroup of the units of this order, and the group of units of any order in $M(2,\bbQ)$ is commensurable with $\SL(2,\bbZ)$.

    To show that $\varrho( \bbZ[\bar\Gamma])$ is an order in $M(2,\bbQ)$. 
    \begin{enumerate}
        \item Let $E=\begin{pmatrix}
      1&1\\
      0&1
    \end{pmatrix},$ $F=\begin{pmatrix}
      1&0\\
      N&1
    \end{pmatrix}$. Then $E$, $F$, $EF$ and $FE\in \varrho( \bbZ[\bar\Gamma])$, so \[\varrho( \bbZ[\bar\Gamma])\otimes_{\bbZ}\bbQ=M(2,\bbQ).\] 
    \item $\varrho( \bbZ[\bar\Gamma])$ is clearly a subring of $M(2,\bbQ)$.
    \item To show that $\varrho( \bbZ[\bar\Gamma])$ is finitely generated $\bbZ$-module, note that $N\varrho( \bbZ[\bar\Gamma])\subset M(2,\bbZ)$, by Lemma~\ref{Structure lemma}.
    \end{enumerate}
\end{proof}

\subsection{Conclusion}
Combining the results from all previous steps, Theorem~\ref{main2} follows.

\section{Trace set of compact hyperbolic surfaces}
\subsection{Proof of Theorem~\ref{main: cocompact 2}}
\begin{proof}[Proof of Theorem~\ref{main: cocompact 2}]
Let $\Sigma_g$ be a closed surface of genus $g\geq 3$ with fundamental group $\Gamma_g$. Let $\{\eta_i\}_{i=1}^{2g}$ be a standard generating set of $\Gamma_g$, and $A$ be the corresponding $(2g-3)$-subgroup.

Recall that the Fricke coordinates for the Teichm\"uller space is a sequence of real numbers $X:=(a_i, c_i, d_i)_{i=1}^{2g-2}$, where $c_i>0$ for all $i$. The embedding corresponding to the sequence $X$ is given by 
\[
\psi_X(\eta_i)=\begin{pmatrix}
    a_i&\frac{a_id_i-1}{c_i}\\
    c_i& d_i
\end{pmatrix}, \quad 1\leq i\leq 2g-2;
\]
\[
\psi_X(\eta_{2g-1})=\begin{pmatrix}
    a&b\\
    c& d
\end{pmatrix}, \quad a+d=b+c>0;
\]
\[
\psi_X(\eta_{2g})=\begin{pmatrix}
    \nu&0\\
    0& \frac{1}{\nu}
\end{pmatrix}, \quad \nu>1.
\]
The numbers $a$, $b$, $c$, $d$, and $\nu$ are uniquely determined up to a sign by the Fricke coordinates and the fundamental relation 
\[
\prod_{i=1}^g[\eta_{2i-1},\eta_{2i}]=e.
\]
For details, see \cite{imayoshi2012introduction}*{page 49} or \cite{hao2022marked}*{Section 8.A}. By abuse of notation, for a point $[d]\in\mathcal{T}_g$, we denote the corresponding embedding from its Fricke coordinates by $\psi_d$. 

Fix $[d]\in\mathcal{T}_g$. Since on the Techm\"uller space the embedding is continuous algebraically, by \cite{mcmullen1999hausdorff}*{Theorem 1.4} and the fact that the Hausdorff dimension of the limit set is equal to the critical exponent for convex-compact Fuchsian groups \cite{patterson1976limit}, the critical exponent of $A$ is a continuous function. Let $V_g^\epsilon(d)$ be an open neighborhood of $[d]$ such that for all points $[d']$ in $V_g^\epsilon(d)$, the critical exponent of $\psi_{d'}(A)$ is greater than $\delta_{\psi_d(A)}-\epsilon$ and Theorem~\ref{critical} holds uniformly for $o$, where $o$ is a fixed base point of ${\bfH}^2$. 

Assume the Fricke coordinates of $[d]$ is given by $(a_i,c_i,d_i)_{i=1}^{2g-2}$. Let $\Theta_{d}$ be the subset of $\mathcal{T}_g$ consisting of elements with Fricke coordinates $(a'_i,c'_i,d'_i)_{i=1}^{2g-2}$ such that $a_i=a_i'$, $b_i=b_i'$ and $c_i=c_i'$ for all $1\leq i\leq 2g-3$. Clearly, for all points in $\Theta_d$, the subgroup $A$ has the same embedding. The construction of the embedding induces a map $f: \Theta_d\to \mathbb{R}$ such that for any point $[d']\in \Theta_d$: 
\[
\psi_{d'}(\eta_{2g})=\begin{pmatrix}
    f([d'])&0\\
    0& \frac{1}{f([d'])}
\end{pmatrix}, \quad f([d'])>1.
\]

Computation shows that $f$ is not constant. Hence, the image of $f$ contains an open neighborhood $I'$ of $f([d])$. 

On the other hand, fix two different elements $a$, $b\in A$, $[d']\in\Theta_d$, and consider the equation:
\[
\mathrm{tr}(\psi_{d'}(\eta_{2g}a))=\mathrm{tr}(\psi_{d'}(\eta_{2g}b)).
\]
Since $\psi_{d'}(a)$ and $\psi_{d'}(b)$ are fixed, this equation in $f([d'])$ has at most two solutions. Let $I$ be the subset of $I'$ that is not a solution of any such equation. For all $[d']\in f^{-1}(I)$, we have:
\[\mathrm{tr}(\psi_{d'}(\eta_{2g}a))\neq \mathrm{tr}(\psi_{d'}(\eta_{2g}b)),\quad
\mathrm{for}\ \mathrm{all}\ a\neq b\in A.\]

Now let $[b']\in V^{\epsilon}_g(d)\cap f^{-1}(I)$. By Theorem~\ref{critical}, there exists a constant $K'>0$ such that for all $R>0$:
\[
\#\setdef{a\in A}{d(o,\psi_{d'}(a)o)\leq R}\geq K'e^{(\delta_{\psi_d(A)}-\epsilon)R}.
\]
Therefore, there exists a constant $K>0$ such that
\[
\#\setdef{a\in A}{d(o,\psi_{d'}(\eta_{2g}a)o)\leq R}\geq Ke^{(\delta_{\psi_d(A)}-\epsilon)R}.
\]
Since
\[
\mathrm{tr}(\psi_{d'}(\eta_{2g}a))\leq 2\cosh(\frac{d(o,\psi_{d'}(\eta_{2g}a)o)}{2}),
\]
the trace set grows at least on the order of $n^{\delta_{\psi_d(A)}-\epsilon}$ for the metric $[d']$.

Now consider a point $[d'']\in V^\epsilon_g(d)\setminus T_{\mathrm{sing}}$. Similarly, there exists a constant $L>0$ such that:
\[
\#\setdef{a\in A}{d(o,\psi_{d''}(\eta_{2g}a)o)\leq R}\geq Le^{(\delta_{\psi_d(A)}-\epsilon)R}
\] 
By the definition of $ T_{\mathrm{sing}}$, $[d'']$ has a minimal marked length pattern, see \cite{hao2022marked}*{Section 8.2}. Therefore, for $[d'']$, we still have: 
\[\mathrm{tr}(\psi_{d''}(\eta_{2g}a))\neq \mathrm{tr}(\psi_{d''}(\eta_{2g}b)),\quad
\mathrm{for}\ \mathrm{all}\ a\neq b\in A.\]

Thus, for $[d'']$, the trace set grows at least on the order of $n^{\delta_{\psi_d(A)}-\epsilon}$.

This completes the proof of Theorem~\ref{main: cocompact 2}.
\end{proof}

\subsection{Proof of Theorem~\ref{main: cocompact}}
\begin{proof}[Proof of Theorem~\ref{main: cocompact}]
Same as in the introduction, let $\Xi_2(g)$ be the set of multicurves on $\Sigma_g$: $\alpha=\cup_{i=1}^s\alpha_i$ suth that all $\alpha_i$ are simple closed geodesics, and $\Sigma_g\setminus \alpha=\Sigma_1\cup \Sigma_2$, where $\Sigma_1$ and $\Sigma_2$ are connected subsurfaces with $|\chi(\Sigma_2)|=2\leq |\chi(\Sigma_1)|.$ Here $|\chi(\Sigma_1)|=2g_1-2+s$, which is the absolute value of the Euler characteristic. For any $[d]\in \mathcal{T}$, denote the length of a multicurve $\alpha$ by $\ell_d(\alpha)=\sum_{i=1}^s \ell_d(\alpha_i)$.

Let $\Delta_g(\epsilon)\subset \mathcal{M}_g$ be the subset such that for any point $(\Sigma_g,[d])$ in $\Delta_g(\epsilon)$, there exists a multicurve $\alpha$ on $\Sigma_g$ such that $\alpha\in \Xi_2(g)$ and $\ell_d(\alpha)\leq \min\{\frac{\pi}{2\kappa}, \frac{\epsilon\pi}{3\kappa}\}$, where $\kappa$ is the constant in Theorem~\ref{buser}.

By \cite{mirzakhani2013growth}*{Theorem 4.9}, we have:
\[
\lim_{g\to \infty} \frac{\mathrm{Vol_{WP}}(\Delta_g(\epsilon))}{\mathrm{Vol_{WP}}(\mathcal{M}_g)}=1.
\]
Since $\mathrm{Vol_{WP}}(T_{\mathrm{sing}})=0$, it is sufficient to show that for all $g\geq 3$, $\Delta_g(\epsilon)\setminus T_{\mathrm{sing}}\subset \bar{V}_g^{\frac{\epsilon}{g}}$.

First, let $[d]\in \Delta_g(\epsilon)$. Note that $\alpha$ cut $\Sigma_g$ into two pieces, $\Sigma_1$ and $\Sigma_2$. Also, $\pi_1(\Sigma_1)=A$ and $\mathrm{Area}(\Sigma_1)=2\pi(2g-4).$ Since $\Sigma_1$ is a compact surface with geodesic boundaries, $\psi_d(A)$ is convex-cocompact, and $\Sigma_1$ is the compact core of ${\bfH}^2/\psi_d(A)$. 

By considering the multicurves $\alpha$, the Cheeger constant satisfies:
\[
h({\bfH}^2/\psi_d(A))\leq\frac{\ell_d(\alpha)}{(4g-8)\pi}.
\]
Now Theorem~\ref{buser} implies that:
\[
\lambda_1({\bfH}^2/\psi_d(A))\leq \frac{\kappa \ell_d(\alpha)}{(4g-8)\pi}=\min\{\frac{1}{8},\frac{\epsilon}{3(4g-8)}\}.
\]
By Theorem~\ref{cheeger}, $\delta_{\psi_d(A)}>\frac{1}{2}$, and 
\[
\lambda_1({\bfH}^2/\psi_d(A))=\delta_{\psi_d(A)}(1-\delta_{\psi_d(A)})>\frac{1}{2}(1-\delta_{\psi_d(A)}).
\]
It follows that: 
\[
\delta_{\psi_d(A)}>\max\{\frac{3}{4}, 1-\frac{2\epsilon}{3(4g-8)}\}.
\]
Now by Theorem~\ref{main: cocompact 2}, $\Delta_g(\epsilon)\setminus T_{\mathrm{sing}}\subset \bar{V}_g^{\frac{4\epsilon}{3(4g-8)}}$.  Since $g\geq 3$, we have $\frac{4\epsilon}{3(4g-8)}\leq \frac{\epsilon}{g}$, and the proof is complete.
\end{proof}

\end{document}